\documentclass[11pt]{amsart}
\usepackage{amscd}
\def\NZQ{\Bbb}               

\def\ZZ{{\NZQ Z}}

\def\frk{\frak}               

\def\Phi{{\frk n}}
\def\Phi{{\frk N}}
\def\opn#1#2{\def#1{\operatorname{#2}}} 
\opn\chara{char} \opn\length{\ell} \opn\pd{pd} \opn\rk{rk}
\opn\projdim{proj\,dim} \opn\injdim{inj\,dim} \opn\rank{rank}
\opn\depth{depth} \opn\grade{grade} \opn\height{height}
\opn\embdim{emb\,dim} \opn\codim{codim}

\opn\Tr{Tr} \opn\bigrank{big\,rank}
\opn\superheight{superheight}\opn\lcm{lcm}
\opn\trdeg{tr\,deg}
\opn\reg{reg} \opn\lreg{lreg} \opn\ini{in} \opn\lpd{lpd}
\opn\size{size}
\opn\div{div} \opn\Div{Div} \opn\cl{cl} \opn\Cl{Cl}
\opn\Spec{Spec} \opn\Supp{Supp} \opn\supp{supp} \opn\Sing{Sing}
\opn\Ass{Ass} \opn\Min{Min}
\opn\Ann{Ann} \opn\Rad{Rad} \opn\Soc{Soc}
\opn\Im{Im} \opn\Ker{Ker} \opn\Coker{Coker} \opn\Am{Am}
\opn\Hom{Hom} \opn\Tor{Tor} \opn\Ext{Ext} \opn\End{End}
\opn\Aut{Aut} \opn\id{id}

\opn\nat{nat}
\opn\pff{pf}
\opn\Pf{Pf} \opn\GL{GL} \opn\SL{SL} \opn\mod{mod} \opn\ord{ord}
\opn\Gin{Gin} \opn\Hilb{Hilb}\opn\sdepth{sdepth}
\opn\aff{aff} \opn\con{conv} \opn\relint{relint} \opn\st{st}
\opn\lk{lk} \opn\cn{cn} \opn\core{core} \opn\vol{vol}
\opn\link{link} \opn\star{star}
\opn\gr{gr}

\def\pot#1#2{#1[\kern-0.28ex[#2]\kern-0.28ex]}

\opn\dirlim{\underrightarrow{\lim}}
\opn\inivlim{\underleftarrow{\lim}}
\let\dirsum=\oplus

\let\iso=\cong
\let\Union=\bigcup

\let\Dirsum=\bigoplus

\let\to=\rightarrow

\def\Implies{\ifmmode\Longrightarrow \else
        \unskip${}\Longrightarrow{}$\ignorespaces\fi}
\def\implies{\ifmmode\Rightarrow \else
        \unskip${}\Rightarrow{}$\ignorespaces\fi}
\def\iff{\ifmmode\Longleftrightarrow \else
        \unskip${}\Longleftrightarrow{}$\ignorespaces\fi}

\let\:=\colon
\newtheorem{Theorem}{Theorem}[section]
\newtheorem{Lemma}[Theorem]{Lemma}
\newtheorem{Corollary}[Theorem]{Corollary}

\let\epsilon\varepsilon
\let\phi=\varphi
\let\kappa=\varkappa
\def\qed{\ifhmode\textqed\fi
      \ifmmode\ifinner\quad\qedsymbol\else\dispqed\fi\fi}
\def\textqed{\unskip\nobreak\penalty50
       \hskip2em\hbox{}\nobreak\hfil\qedsymbol
       \parfillskip=0pt \finalhyphendemerits=0}
\def\dispqed{\rlap{\qquad\qedsymbol}}

\opn\dis{dis}
\def\pnt{{\raise0.5mm\hbox{\large\bf.}}}

\opn\Lex{Lex}

\begin{document}

\title{Stanley decompositions, pretty clean filtrations and reductions modulo regular elements}

\author{Asia Rauf}
\thanks{}
\address{Asia Rauf,
School of Mathematical Sciences, GCU, Lahore}
\email{asia.rauf@gmail.com}

\begin{abstract}
We study the behavior of Stanley decompositions and of pretty clean
filtrations under reduction modulo a regular element.
\end{abstract}

\maketitle

\section*{Introduction}
Let $K$ be a field and $S=K[x_1,x_2,\ldots,x_n]$ be a polynomial
ring in $n$ variables over the field $K$. Let $I\subset S$ be a
monomial ideal, and let $u\in S$ be a monomial such that $u$ is
regular on $S/I$. The purpose of this paper is to investigate how
the Stanley depth and the property of $S/I$ to be pretty clean
behaves  when we pass from $S/I$ to $S/(I,u)$, and vice versa.

We denote by $I^c\subset S$ the $K$-linear subspace
of $S$ generated by all monomials which do not belong to $I$. Then
$S=I\bigoplus I^c$  and
 $S/I \cong I^c$ as $K$-linear spaces.

If $u\in S$ is a monomial and $Z\subset \{x_1,\ldots,x_n\}$, the
$K$-subspace $uK[Z]$ whose basis consists of all monomials $uv$,
with $v\in K[Z]$, is called a {\em Stanley space} of dimension
$|Z|$. A decomposition $\mathcal{D}$ of $I^c$ as a finite direct sum
of Stanley spaces is called a {\em Stanley decomposition} of $S/I$.
The minimal dimension of a Stanley space in $\mathcal{D}$ is called
the {\em Stanley depth} of ${\mathcal D}$ and it is denoted by
$\sdepth({\mathcal D})$. We set
\[
\sdepth(S/I):=\max\{\sdepth({\mathcal D}):\; {\mathcal D}\; \text{is a Stanley
decomposition of}\;  S/I \}
\]
 and  call this number the {\em Stanley depth} of $S/I$.

Stanley \cite[Conjecture 5.1]{St1} made a conjecture on general Stanley decompositions of $\ZZ^n$-graded modules. In the special case that the $\ZZ^n$-graded module is $S/I$, where $I$ is a monomial ideal,  the conjecture says that
 $\sdepth(S/I)\geq \depth(S/I)$. A monomial ideal $I$ is called {\em Stanley ideal} if
it satisfy  Stanley's conjecture.

A basic fact in commutative algebra says   $\depth S/(I,f)=\depth
S/I-1$ for any homogeneous element of positive degree  $f\in S$
which is regular on $S/I$. In this paper we  show that a
corresponding statement holds for the Stanley depth. In fact, we
show in Theorem~\ref{main} that $\sdepth (S/(I,u))=\sdepth S/I-1$
for any monomial $u\in S$  which is regular on $S/I$.

Special Stanley decompositions arise from prime filtrations. Let
$$\mathcal F: I=I_0 \subset I_1 \subset \ldots \subset I_r = S$$ be
a prime filtration of $S/I$, i.e. $I_j/I_{j-1}\cong S/P_j$ for any
$j=1,\ldots,r$, where $P_j\subset S$ are prime ideals. The support
of $\mathcal F$, is the set $\Supp(\mathcal F)=\{P_1,\ldots,P_r\}$.
It is well known that $\Ass(S/I)\subset \Supp (\mathcal F)$. Recall
that the prime filtration $\mathcal F$ is called {\em pretty clean},
if for all $i<j$ with $P_i\subset P_j$ it follows that $P_i=P_j$. If
$S/I$ has a pretty clean filtration then $S/I$ is called pretty
clean, see \cite[Definition 3.3]{HP}. For the pretty clean
filtration, $\Supp(\mathcal F)=\Ass(S/I)$, see \cite[Corollary
3.6]{HP}. This condition implies, by \cite[Proposition 2.2]{A}, that
$I$ is a Stanley ideal. The prime filtration $\mathcal F$ is {\em
clean} if $\Supp(\mathcal F)=\Min(S/I)$, where $\Min(S/I)$ is the
set of minimal prime ideals of $S/I$. Note that any clean filtration
is pretty clean. If $S/I$ has a  clean filtration, that $S/I$ is
called clean.

The main result (Theorem \ref{lahore}) of the second section is that
if $I\subset S$ is a monomial ideal and $u\in S$ is a monomial which
is regular on $S/I$, then $S/I$ has a pretty clean filtration if and
only if $S/(I,u)$ has a pretty clean filtration. This result implies
that an ideal generated by a regular sequence of monomials is pretty
clean. This fact was first proved in \cite[Proposition 1.2]{HSY} by
a different  method.

I want to thank Professor J\"urgen Herzog for his advice during the preparation of the paper.

\section{Stanley decompositions and regular elements}

The aim of this section is to show that the Stanley depth behaves
like the ordinary depth with respect to reduction modulo regular elements. Indeed we
have the following result:

\begin{Theorem}
\label{main}
Let $I\subset S$ be a monomial ideal of $S=K[x_1,\ldots,x_n]$ and
$u\in S$ be a monomial regular on $S/I$. Then
$\sdepth(S/(I,u))=\sdepth(S/I)-1.$ In particular, $I$ is a Stanley
ideal if and only if $(I,u)$ is a Stanley ideal.
\end{Theorem}

We first prove a special case of the theorem:

 \begin{Lemma}
\label{newvariable} Let $m<n$ and  $J\subset S'=K[x_1,\ldots,x_m]$
be a monomial  ideal.  Then for the monomial ideal $I=JS$ and for
any $x_k$ with $m<k\leq n$ we have
$$\sdepth(S/(I,x_k)) = \sdepth(S/I) - 1.$$
\end{Lemma}

\begin{proof}
Let $T=S'[x_{m+1},\ldots, x_{k-1},x_{k+1}\ldots,x_n]$  and $L\subset T$ be
the monomial ideal such that $L=JT$. Then we have $S/(I,x_k)=T/L$. Let
\[
 \mathcal{D}:\; T/L=\Dirsum_{i=1}^ru_iK[Z_i]
 \]
be a Stanley decomposition of $T/L$ such that $\sdepth
\mathcal{D}=\sdepth T/L$. Then
\[
\mathcal{D}_1: S/I =
 (T/L)[x_k]=\Dirsum_{i=1}^ru_iK[Z_i][x_k]=\Dirsum_{i=1}^ru_iK[Z_i,x_k].
 \]
 is a Stanley decomposition of $S/I$.
 It follows that
 \[
 \sdepth {\mathcal{D}_1} = \sdepth {\mathcal{D}}+1= \sdepth T/L+1
  \]   and
 \[ \sdepth {\mathcal{D}_1 } \leq  \sdepth S/I.
 \]
 Hence
 \[
\sdepth T/L+1 \leq \sdepth S/I. \] In order to prove the opposite
inequality we consider a Stanley decomposition
\[
\mathcal{D}_2:\; S/I=\Dirsum_{i=1}^sv_iK[W_i]
\]
of $ S/I$ with $\sdepth \mathcal{D}_2=\sdepth S/I$.

Let $\mathcal I=\{i\in [s] :\;  v_iK[W_i] \cap T\neq \{0\}\}$.
  We claim that

\begin{eqnarray}
\label{sect1} \mathcal{D}_3 :\; T/L=L^c=\Dirsum_{i\in \mathcal
I}v_iK[W_i] \cap T.
\end{eqnarray}
and $\Dirsum_{i\in \mathcal I}v_iK[W_i] \cap T$ is a  direct sum
decomposition of $T/L$.

  In order to prove (\ref{sect1}), choose a monomial $v\in L^c$. We want to show that there exists $i\in \mathcal I$
such that $v\in v_iK[W_i] \cap T$. Suppose on the contrary that
$v\not\in v_iK[W_i] \cap T$ for all $i\in \mathcal I$. Since $v\in T$, it implies
that $v\not\in v_iK[W_i]$, for all $i$. Hence we have $v\in I=JS$.
Since $v\in T$ and $L=JT$, it follows that $v\in L$, a
contradiction. Conversely, choose a monomial $w\in v_iK[W_i] \cap
T$.  This implies that
$w\not\in I=JS$ and since $L=JT\subset JS=I$, we see that  $w\in L^c$.

 Now we will show that $\mathcal{D}_3$ is a Stanley decomposition. Indeed, we have
\[
 v_i K[W_i]\cap T = \begin{cases}
v_i K[W_i\setminus \{x_k\}],& \text {if}\; x_k\text{ does not divide } v_i\\
{0} ,& \text{otherwise.} \end{cases}
\]

Comparing the Stanley decomposition ${\mathcal D}_2$ of $S/I$ with the  Stanley decomposition ${\mathcal D}_3$ of $T/L$ we see that
$\sdepth(\mathcal{D}_2)\leq \sdepth(\mathcal{D}_3)+1$. Hence
\[
\sdepth S/I=\sdepth \mathcal{D}_2\leq \sdepth(\mathcal{D}_3)+1\leq
\sdepth T/L + 1.
\]
\end{proof}

For the proof of  Theorem \ref{main} we also need the following simple fact:

\begin{Lemma}
\label{easy}
Let \[ I=I_0\subset I_1\subset \ldots \subset I_r=S
 \]
be an ascending chain of monomial ideals of $S$ such that each
$I_j/I_{j-1}$ is a cyclic module, and hence $I_j/I_{j-1}\iso
S/L_j(-a_j)$ for some monomial ideal $L_j$ and some $a_j\in \ZZ^n$.
Then
\[
 \sdepth(S/I)\geq \min \{ \sdepth(S/L_j) :\; j\in \{1,\ldots,r\} \}
\]
\end{Lemma}

\begin{proof} We have the following decomposition of $S/I$ as a $K$-vector space:
\[
S/I=I_1/I_0\dirsum I_2/I_1\dirsum \cdots \dirsum S/I_{r-1}.
\]
Since each $I_j/I_{j-1}\iso S/L_j(-a_j)$ we get the isomorphism
\begin{eqnarray}
\label{shift}
S/I\iso S/L_1(-a_1)\dirsum S/L_2(-a_2)\dirsum \cdots  \dirsum S/L_r(-a_r).
\end{eqnarray}
For each $j$ let ${\mathcal D}_j:  S/L_j= \Dirsum_{k=1}^{r_j}u_{jk}K[Z_{jk}]$ be a Stanley decomposition of $S/L_j$ such that $\sdepth {\mathcal D}_j=\sdepth S/L_j$. Then by the isomorphism (\ref{shift}) we obtain the following Stanley decomposition
\[
S/I= \Dirsum_{j=1}^r\Dirsum_{k=1}^{r_j}u_ju_{jk}K[Z_{jk}],
\]
of $S/I$, where $u_j=x^{a_j}$ for $j=1,\ldots,r$. From this Stanley decomposition of $S/I$ the desired inequality follows.
\end{proof}

\begin{proof}[Proof of Theorem \ref{main}]
Without loss of generality we may assume that $I=JS$ where $J\subset S'=K[x_1,\ldots,x_m]$
and that  $u=x_{m+1} ^{a_1}\ldots
x_n ^{a_{n-m}}$. We consider an ascending chain of ideals of $S$ between $(I,u)$ and $S$ where two successive members of the chain are of the form
\[
(I,x_{m+1} ^{b_1}\cdots x_{k} ^{b_k} \cdots x_n
^{b_{n-m}})\subset(I,x_{m+1} ^{b_1}\cdots x_{k}^{b_{k}-1} \cdots
x_n ^{b_{n-m}})
\]
and where $b_i\leq a_i$ for all $i=1,\ldots, n-m$.

Observe that
\[
(I,x_{m+1} ^{b_1}\cdots x_{k} ^{b_{k-1}} \cdots x_n ^{b_{n-m}})/
(I,x_{m+1} ^{b_1}\cdots x_{k} ^{b_k} \cdots x_n^{b_{n-m}}) \simeq
S/(I,x_k).
\]
Therefore Lemma \ref{newvariable} and Lemma \ref{easy} imply that
\[
\sdepth(S/(I,u))\geq  \sdepth(S/(I,x_k)) = \sdepth(S/I)-1.
\]
In order to prove other inequality, we choose a Stanley decomposition
$$\mathcal{D}': (I,u)^c= \bigoplus_{i=1}^r u_i K[Z'_i]$$  of
$S/(I,u)$ with $\sdepth(\mathcal{D}')=\sdepth(S/(I,u))$. We obtain a
direct sum of $K$-vector subspaces $\bigoplus_{i=1}^r u_i
K[Z'_i]\cap S'$ of $S'$. We observe that
\[
J^c= \bigoplus_{i=1}^r u_i K[Z'_i]\cap S'
\]
and that $\bigoplus_i u_i K[Z'_i]\cap S'$ is  a Stanley
decomposition of $S'/J$, where the sum is taken over those $i\in
\{1,\ldots, r\}$ for which $u_i K[Z'_i]\cap S'\neq \{0\}$, cf.\  proof of Lemma~\ref{newvariable}.

We have
\[
 u_i K[Z'_i]\cap S' = \begin{cases}
u_i K[Z'_i\cap \{x_1,\ldots,x_m \}],& \text {if}\; \supp(u_i)\subset \{x_1,\ldots,x_m\}\\
{0} ,& \text{otherwise.} \end{cases}
\]
Hence if we set  $\Lambda=\{i : \; \supp(u_i)\subset
\{x_1,\ldots,x_m\} \}$, then
\[
\mathcal{D}: S/I = \bigoplus_{i\in \Lambda} u_i K[Z_i]
\]
is a Stanley decomposition of $S/I$, where $Z_i:=\{Z'_i\cap
\{x_1,\ldots,x_m \}\}\cup \{x_{m+1},\ldots,x_n\}$.

We claim that $|Z_i|>|Z'_{i}|$. Indeed, otherwise
$\{x_{m+1},\ldots,x_n\} \subset Z'_i$, contradicting the fact that
$(u)\cap u_i K[Z'_i]=\{0\}$. Therefore, $\sdepth(\mathcal{D})\geq
\sdepth(\mathcal{D}')+ 1$.

Hence or final conclusion is that
\[
\sdepth(S/(I,u))= \sdepth(S/I)-1.
\]
\end{proof}

As an immediate consequence of our theorem we obtain the following
result first proved in \cite[Proposition 1.2]{HSY}.

\begin{Corollary}
\label{regular}
Let $I$ be a monomial ideal generated by regular sequence of monomials. Then $I$ is a Stanley ideal.
\end{Corollary}

\section{Pretty clean filtrations and regular elements}

\begin{Theorem}
\label{lahore} Let $S=K[x_1,x_2,...,x_n]$  be a polynomial ring and
$I\subset S$ be a  monomial ideal and  $u$  a monomial in $S$ such
that $u$ is regular on $S/I$. Then  $S/I$ is pretty clean if and
only if $S/(I,u)$ is pretty clean.
\end{Theorem}

\begin{proof} Suppose $S/I$ is pretty clean and let
\[
\mathcal{F}:\; I=I_0\subset
                  I_1\subset\ldots \subset I_r=S
                  \]
be a pretty clean  filtration of $S/I$ with   $I_j/I_{j-1}\cong S/P_j$ for
$j=1,2,...,r$. It is known from \cite[Corollary 3.6]{HP} that
$\Ass(S/I)=\{P_1,\ldots, P_r\}$.

We have  $I_j=(I_{j-1},z_j)$ where $z_j$ is a monomial in $S$. The
prime filtration $\mathcal{F}$ induces the following filtration
\[
\mathcal{G}:\; (I,u)\subset
                  (I_1,u)\subset \ldots\subset (I_r,u)=S,
                  \]
where
\[
(I_j,u)/(I_{j-1},u)=((I_{j-1},u),z_j)/(I_{j-1},u)\cong
S/(I_{j-1},u):z_j.
\]
Since $u$ is regular on $S/I$, it follows that  $u$ is regular on
$S/I_j$ for all $j$.  Indeed, since $S/I$ is pretty clean it follows
that $S/I_j$ is pretty clean. Hence
$\Ass(S/I_j)=\{P_{j+1},\ldots,P_r\}$ which is contained in
$\Ass(S/I)$. Since $\gcd(u,z_j)=1$ it follows that
$$(I_{j-1},u):z_j=((I_{j-1}:z_j),u)=(P_j,u).$$
Hence
\[
(I_j,u)/(I_{j-1},u)\cong S/(P_j,u).
\]
Suppose, without loss of generality, that
\[
P_j=(x_1,...,x_t)\quad \text{and}\quad  u=\prod_{i=t+1}^n
{x_i}^{a_i}.
\]
Then $S/(P_j,u)\cong K[x_{t+1},...,x_n]/(u)K[x_{t+1},...,x_n]$,
which is clean by \cite{A}. Hence we see that $(I_j,u)/(I_{j-1},u)$
is clean and $$\Ass((I_j,u)/(I_{j-1},u))=\{(P_j,x_i) : \; x_i\mid
u\}.$$

Therefore our filtration  $\mathcal{G}$  can be refined as follows
\[
(I_{j-1},u)=I_{{j-1},0}\subset I_{{j-1},1}\subset \ldots\subset
I_{{j-1},s_j}=(I_j,u)
\]
where
\[ I_{{j-1},k}/I_{{j-1},k-1}\cong
S/P_{j-1,k}
\]
with $P_{j-1,k}\in\{(P_j,x_{i}): x_{i}\mid u\}$.

In the refined filtration of  $\mathcal{G}$ if we have
$I_{j,k}\subset I_{i,l}$, then either $j<i$ or $j=i$ and $k<l$.
Suppose $j<i$ and $P_{j,k}\subset P_{i,l}$. We have
$P_{j,k}=(P_{j+1},x_r)$ for some $r$ and $P_{i,l}=(P_{i+1},x_s)$ for
some $s$. Since $u\not\in \Union_{P\in \Ass(S/I)}P$ it follows that
$x_s\not\in P_{j+1}$. Therefore, $P_{j+1}\subseteq
P_{i+1}$. However, since ${\mathcal F}$ is a pretty clean filtration
it follows that $P_{j+1}=P_{i+1}$, and hence $P_{j,k}= P_{i,l}$.

Next suppose that $i=j$ and $k<l$ and suppose that $P_{i,k}\subseteq P_{i,l}$. Since $\height P_{i,k}=\height P_{i,l}$
we conclude that $P_{j,k}= P_{i,l}$, also in this case. Thus we have shown that the refinement of $\mathcal G$ is a pretty clean filtration of  $S/(I,u)$, and hence $S/(I,u)$ is pretty clean.

Conversely, suppose that $S/(I,u)$ is pretty clean. Since $u$ is regular on $S/I$, we
may suppose that $I=JS$ where $J\subset S'=K[x_1,\ldots,x_m]$ for
$m<n$ and $\supp(u)\subset \{x_{m+1},\ldots,x_n \}$. Since $S/(I,u)$ is
pretty clean there exist a pretty clean filtration
\[
\mathcal{M}:\; (I,u)=I'_0\subset
                  I'_1\subset\ldots \subset I'_r=S
                  \]
such that $I'_j/I'_{j-1}\cong S/P_j$ where $P_j\in \Ass(S/(I,u))$.
Recall that
\[
\Ass(S/(I,u))=\{(P',x_k) \: P'\in \Ass(S'/J) \text{ and }  x_k \mid
u\}.
\]
By taking the intersection of above filtration
$\mathcal{M}$ with $S'$, we get the  filtration
\[
 \mathcal{N}:\; J_0\subseteq
                  J_1\subseteq \ldots\subseteq J_r=S'
                  \]
of $S'/J_0$ where  $J_j=I'_j\cap S'$ for $j=0,\ldots,r$. We claim
that $J_0=J$. Let $I$ be generated by the monomials
$u_1,\ldots,u_l$. Since $I=JS$ with $J\subset S'$ it follows that
$u_i\in S'$ for all $i$. Choose a monomial $v\in J_0=(I,u)\cap S'$.
Then either $v=eu_i$ where $e\in S'$, or  $v=fu$ where $f\in S'$.
The second case cannot happen since $v\in S'$. This shows that
$J_0\subset J$. The other inclusion is obvious.

 Take an ideal $I'_j\in  \mathcal{M}$. Then   $I'_j=(I'_{j-1},w_j)$
 where $w_j\in S$ and $(I'_{j-1}:w_j)=(P',x_k)$ for some $P'\in \Ass(S'/J)$ and some $x_k$ such that $x_k\mid u$. Then we have $I'_{j-1}\cap S'=I'_j\cap
 S'$ if and only if $w_j\not\in S'$.

Let  $\{r_0,\ldots,r_k\}$ be the subset of $[r]$ for which we have
$J_{r_{i}}$ is properly contained in $J_{r_{i}+1}$ in the filtration
$\mathcal N$. Set $L_i=J_{r_{i}}$ for $i=0,\ldots, k$ and
$L_{k+1}=S'$. Then we obtain the  filtration
\[
\mathcal{L}:\; J= L_0\subset
                  L_1\subset\ldots \subset L_{k+1}=S'.
\]
We note that $L_i=(J,w_{{r_0}+1},w_{{r_1}+1},\ldots, w_{r_{i-1}+1})$
for $i=0,\ldots,k+1$ with $w_{{r_i}+1}\in S'$ for all $i$.

Since  $L_i=(L_{i-1},w_{r_{i-1}+1})$, we have that  $L_i/L_{i-1}\cong S'/(L_{i-1}:_{S'} w_{r_{i-1}+1})$ and also we have that  $L_i=I'_{r_{i}}\cap S'$. So
  $(L_{i-1}:_{S'} w_{r_{i-1}+1})=(I'_{r_{i-1}} \cap S':_{S'} w_{r_{i-1}+1})$.

  We claim that
 $(I'_{{r_{i-1}}}\cap S':_{S'} w_{r_{i-1}+1})=(I'_{{r_{i-1}}}:_{S} w_{r_{i-1}+1})\cap
 S'$. In fact, the inclusion $(I'_{{r_{i-1}}}\cap S':_{S'} w_{r_{i-1}+1})\subset(I'_{{r_{i-1}}}:_{S} w_{r_{i-1}+1})\cap
 S'$ is obvious. In order to prove the other inclusion  we choose a monomial $v\in (I'_{{r_{i-1}}}:_{S} w_{r_{i-1}+1})\cap
 S'$.
Then we have that  $v\in (I'_{{r_{i-1}}}:_{S} w_{r_{i-1}+1})$ and
$v\in S'$. Hence $vw_{r_{i-1}+1}\in I'_{{r_{i-1}}}$ and
$vw_{r_{i-1}+1}\in S'$, since $w_{r_{i-1}+1}\in S'$. Therefore
$vw_{r_{i-1}+1}\in I'_{{r_{i-1}}}\cap S'$ which implies that $v\in
(I'_{{r_{i-1}}}\cap S':_{S'} w_{r_{i-1}+1})$, as desired.

Now we see that
\begin{eqnarray*}
(L_{i-1}:_{S'} w_{r_{i-1}+1})&=&(I'_{{r_{i-1}}}\cap S':_{S'}
w_{r_{i-1}+1})\\
&=&(I'_{{r_{i-1}}}:_{S} w_{r_{i-1}+1})\cap S'
=(P',x_k)\cap
S'=P',
\end{eqnarray*}
where $(P',x_k)\in \Ass(S/(I,u))$.

This shows that $\mathcal L$ is a prime filtration with the property
that the prime ideals in  $\Supp({\mathcal L})$ form a subsequence
of $P_1,\ldots, P_r$. Therefore, since $\mathcal{M}$  is a pretty
clean filtration, the filtration $\mathcal L$ is pretty clean as
well.  From this fact we will deduce that $S/I$ is pretty clean.
This then will complete the proof of the theorem.

Indeed, our filtration $\mathcal L$ induce the filtration
\[
{\mathcal K}: I=JS=L_0S\subset L_1S\subset\cdots \subset L_{k+1}S=S.
\]
with $L_iS/L_{i-1}S\iso S/P'S$ where $L_i/L_{i-1}\iso S'/P'$ for
$i=1,\ldots,k+1$. This holds because the extension $S'\to S$ is
flat. Now, since $\mathcal L$ is a pretty clean filtration of
$S'/J$,  it is obvious that $\mathcal K$ is  a pretty clean
filtration of $S/I$.
\end{proof}

As an immediate consequence we obtain the following result from
\cite[Proposition 1.2]{HSY}.

\begin{Corollary}
\label{essen} Let $u_1,\ldots,u_k$ be a regular sequence  in the
polynomial ring $S$. Then $S/(u_1,\ldots,u_k)$ is pretty clean.
\end{Corollary}

\begin{proof}
We use induction on $k$. For $k=1$ the assertion follows from
Theorem \ref{lahore} applied to $I=(0)$, or from \cite{A}. By
induction hypothesis we may now assume that $S/(u_1,\ldots,u_{k-1})$
is pretty clean.  Since $u_k$ is regular on $S/(u_1,\ldots,u_{k-1})$
it follows again from Theorem \ref{lahore} that $S/(u_1,\ldots,u_k)$
is pretty clean.
\end{proof}


\newpage

\begin{thebibliography}{99}
\bibitem[1]{HSY} J.\ Herzog, A.\ Soleyman Jahan and S.\ Yassemi, Stanley decompositions and partitionable simplicial complexes,   math.AC/0612848, to appear in J.\ Alg.\ Comb.\ 2007.
\bibitem[2]{HP} J.\ Herzog and D.\ Popescu, Finite Filtrations of Modules and Shellable
Multicomplexes, manuscripta math.\ {\bf 121} (2006), 385--410.
\bibitem[3]{A} A. Soleyman Jahan, Prime filtrations of monomial ideals and polarizations,  J.\ Alg.\  {\bf 312}(2) (2007), 1011--1032.
\bibitem[4]{St1} R.\ P.\ Stanley, Linear Diophantine Equations  and Local Cohomology, Invetn.\ math.\ {\bf 68} (1982), 175--193.
\end{thebibliography}
\end{document}